# Robust adaptive droop control for DC microgrids

Tuyen V. Vu [1*], Dallas Perkins [1], Fernand Diaz [1], David Gonsoulin [1], Chris S. Edrington [1], Touria El-Mezyani [2]

[1] Center for Advanced Power Systems, Florida State University, Tallahassee, FL, USA

[2] Department of Electrical Engineering, University of West Florida, Pensacola, FL, USA

ABSTRACT

There are tradeoffs between current sharing among distributed resources and DC bus voltage stability when conventional droop control is used in DC microgrids. As current sharing approaches the setpoint, bus voltage deviation increases. Previous studies have suggested using secondary control utilizing linear controllers to overcome drawbacks of droop control. However, linear control design depends on an accurate model of the system. The derivation of such a model is challenging because the noise and disturbances caused by the coupling between sources, loads, and switches in microgrids are under-represented. This under-representation makes linear modeling and control insufficient. Hence, in this paper, we propose a robust adaptive control to adjust droop characteristics to satisfy both current sharing and bus voltage stability. First, the time-varying models of DC microgrids are derived. Second, the improvements for the adaptive control method are presented. Third, the application of the enhanced adaptive method to DC microgrids is presented to satisfy the system objective. Fourth, simulation and experimental results on a microgrid show that the adaptive method precisely shares current between two distributed resources and maintains the nominal bus voltage. Last, the comparative study validates the effectiveness of the proposed method over the conventional method.

## 1. Introduction

### 1.1. Literature Review

Microgrids play a very important role in the future smart grid as they provide clean, reliable, and affordable energy solutions to both utilities and electrical users [1]-[5]. When comparing the two types of microgrids , DC type have competitive advantages over AC type regarding their low distribution loss and simple control algorithms [6]-[9].

The traditional control methodology in DC microgrids is voltage droop control. The control objective of these systems, to share the current among distributed resources (DR), is achieved through reducing the output voltage of the resources following predefined droop characteristic lines. Various types of droop control can be found in the literature [10]. In [10], the authors review five types of droop control for current sharing among DR. The droop control method is widely utilized for current sharing among DR because of its simplicity in implementation; however, the limitation of the method is that the enhancement in current sharing increases the deviation of DC bus voltage from its nominal value [11]-[18]. On one hand, increased droop resistances result in an improved accuracy in current sharing but increased bus voltage deviation. On the other hand, decreased droop resistances result in reduced bus voltage deviation but inaccurate current sharing.

There have been efforts to improve both current sharing and bus voltage stability. Latest efforts based on secondary control systems are the utilization of linear control for (a) bus voltage compensation and (b) current sharing compensation. However, the system models used for the linear control design using stability criteria has yet to be fully developed [11]-[25]. Another control trend is to find the desired droop values, which satisfy the system's objective, based on rules for changing droop resistances or linear controllers to change the droop resistances [26]-[30]. It is apparent that the rule-based methods in [26]-[30] do not present an automatic and/or optimal approach in control systems. Thus, the compensation in such methodologies is not accurate. The improvement using linear controllers to change droop characteristics in [30], [31] requires an accurate model. Such an accurate model is not achievable in DC microgrids because (a) there are electrical couplings among DR, load devices, and distribution systems; and (b) there are load disturbances, sensor noise, and plug-n-play phenomenon in DR.

Thus, this paper proposes a control methodology that does not require an accurate model and but that is robust to the system's uncertainties. According to [32], among control methodologies, the adaptive control does not require in-depth knowledge about the system for control design. The control parameters will self-adapt to achieve the optimal performance. Therefore, the adaptive control has been selected for investigation and analysis in this paper. Adaptive control has evolved over the last four decades [33]-[35]. The developments seemed to slow down in the 1990s with the limitations of conventional adaptive control in the sense that faster adaptation rates produce more oscillatory state responses. However, the interest in adaptive control has been renewed in the past decade with modifications in adaptive control structure for transient improvements. These modifications for advancements in adaptive control can be found in [36]-[43]. The first improvement results in fast adaptation

NOMENCLATURE

| | |
|---|---|
| $a$ | Unknown system parameter. |
| $a_m$ | Reference model parameter. |
| $a_{Ii}$ | Unknown current model parameter. |
| $a_{Imi}$ | Reference current model control parameter. |
| $a_{Vi}$ | Unknown voltage model parameter. |
| $a_{Vmi}$ | Reference voltage model control parameter. |
| $b$ | Unknown system's parameter. |
| $b_{Ii}$ | Unknown current model parameter. |
| $b_{Imi}$ | Reference current model control parameter. |
| $b_{Vi}$ | Unknown voltage model parameter. |
| $b_{Vmi}$ | Reference voltage model control parameter. |
| $b_m$ | Reference model parameter. |
| $e$ | Error signal. |
| $\hat{e}$ | Estimate of error signal $e$. |
| $e_{Ii}$ | Current error signal. |
| $\hat{e}_{Ii}$ | Estimate of error signal $e_{Ii}$. |
| $e_{Vi}$ | Voltage error signal. |
| $\hat{e}_{Vi}$ | Estimate of error signal $e_{Vi}$. |
| $G_{Ci}(s$ | Converter voltage transfer function. |
| $G_{Ii}(s)$ | Current to droop transfer function. |
| $G_{Vi}(s)$ | Bus voltage to droop transfer function. |
| $G_{Imi}(s)$ | Reference current control model. |
| $G_{Vmi}(s)$ | Reference voltage control model. |
| $g(\bar{\theta})$ | Boundary of convex set $M_{\bar{\theta}}$. |
| $g(b)$ | Boundary of convex set $M_b$. |
| $I$ | Identity matrix. |
| $i_i$ | Distributed resource current. |
| $i_{mi}$ | Output of the current reference control model. |
| $i_L$ | Load current. |
| $\hat{\imath}_i$ | Small signal of distributed resource current. |
| $i_{i,pu}$ | Per-unit current control feedback. |
| $i_{i,pu}^{ref}$ | Per-unit current control reference. |
| $l$ | Error feedback gain. |
| $l_{Ii}$ | Error feedback gain of current reference model. |
| $l_{Vi}$ | Error feedback gain of voltage reference model. |
| $M_{\bar{\theta}}$ | Convex set of $\bar{\theta}$. |
| $M_b$ | Convex set of $b$. |
| $m_{Ii}$ | Normalized term of current control. |
| $m_{Vi}$ | Normalized term of bus voltage control. |
| NPC | Neutral point clamped converter. |
| $R_{di}$ | Droop resistance. |
| $R_{Ii}$ | Droop change regarding current control. |
| $R_{Vi}$ | Droop change regarding bus voltage control. |
| $\hat{r}_{di}$ | Small signal of $R_{di}$. |
| $V_i$ | Distributed resource output voltage. |
| $V_{bi}$ | Terminal bus voltage. |
| $V_i^{ref}$ | Reference voltage for distributed resource. |
| $\bar{v}_{i,pu}$ | DC microgrid average bus voltage. |
| $v_{pu}^{ref}$ | Average reference bus voltage. |
| $v_{mi}$ | Output of the bus voltage reference control model. |
| $\hat{v}_i$ | Small signal of resource output voltage. |
| $\hat{v}_{bi}$ | Small signal of terminal bus voltage. |
| $u$ | Control signal. |
| $u_n$ | Normalized signal of $u$. |
| $x$ | System's state variable. |
| $x_m$ | Reference model state variable. |
| $\bar{\theta}$ | Adaptive control parameter vector. |
| $\bar{\theta}_{Ii}$ | Adaptive current control parameter vector. |
| $\bar{\theta}_{Vi}$ | Adaptive voltage control parameter vector. |
| $\phi$ | Signal vector. |
| $\phi_n$ | Normalized signal of $\phi$. |
| $\phi_{Ini}$ | Normalized signal of current control. |
| $\phi_{Vni}$ | Normalized signal of bus voltage control. |
| $\gamma$ | Adaptation gain. |
| $\gamma_k$ | Scheduled adaptation gain. |
| $\gamma_{Iik}$ | Scheduled adaptation gain of current control. |
| $\gamma_{Vik}$ | Scheduled adaptation gain of bus voltage control. |
| $\varepsilon$ | Modeling error. |
| $\varepsilon_{Vi}$ | Modeling error of current control. |
| $\varepsilon_{Ii}$ | Modeling error of bus voltage control. |

rate and improved transient response [37]. Nevertheless, the assumption of initial state errors, which are assumed to be zero, are always violated. Moreover, the use of $\mathcal{L}_\infty$ norms to determine the boundedness in the transient performance is not practical. Other improvements based on the closed-loop reference model (CRM) are proposed in [38] and further investigated in [39]-[43], which prove the stability, robustness, and transient performance of using $\mathcal{L}_2$ norms analysis. However, the recent advancements do not guarantee that all the signals in the system are bounded, which may cause the system to be unstable. Thus, additional improvements to make all the signals in the systems bounded are presented in this paper. Moreover, using a fixed adaptation gain in the method for various operation modes, which is popular in DC microgrids, results in substantial differences in the output response during the transient and may lead the system into instability. Hence, a modification in the adaptation gain is made to ensure that the system is stable at various operation modes.

### 1.2. *Contributions of the Paper*

As mentioned, previous studies have not addressed the controls design for DC microgrids because of the inability to derive accurate system models. Therefore, in this paper, we first present the derivation of the time-varying models for DC microgrids. Secondly, the paper presents the analysis and enhancement of the closed-loop model reference adaptive control utilizing (a) normalization technique and (b) projection algorithm for robustness of the control systems. Moreover, an adaptation gain-scheduling technique to improve the transient response in various voltage or current operating points is proposed. Third, the developed control method is utilized in a distributed secondary control layer (DSC) to iteratively find the desired but unknown droop characteristics of each DR in DC microgrids for the first time. Last, the simulation and experiments are conducted in the secondary control layer of a DC microgrid with communication among distributed adaptive controllers to validate the significance of the proposed adaptive control in ensuring the stability and robustness of the system.

The organization of the paper is as follows: Section II introduces the proposed control methodology, which includes (a) the derivation of a small-signal model based on conventional droop method, (b) the improvement of an adaptive control method, and (c) adaptive control application for adaptive droop in DC microgrids. The simulation and experiments of a DC microgrid are presented in Section III. The results are analyzed and discussed in Section IV to verify the effectiveness of the method. Section V elaborates the advantages of the proposed method over the linear control method. Section VI concludes the achievements of the paper.

## 2. Proposed Adaptive Droop Methodology

This section is formatted as follows: First, linear time-varying (LTV) models for DC microgrids based on droop control are formulated. Second, a robust adaptive control algorithm is developed and analyzed. Finally, the developed adaptive control and the LTV models are utilized for a proposed adaptive droop control applied to control the current sharing among DR and to stabilize the bus voltage in DC microgrids.

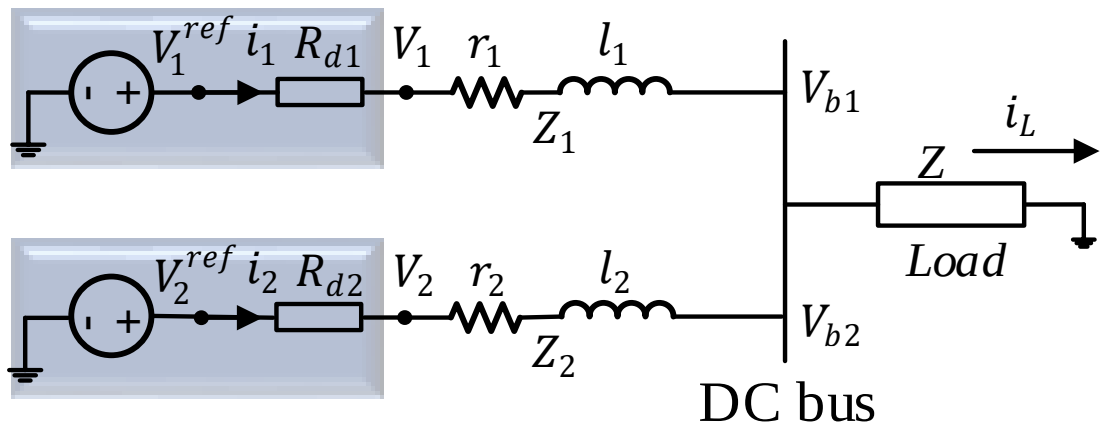

Fig. 1. Example of parallel DR with droop control.

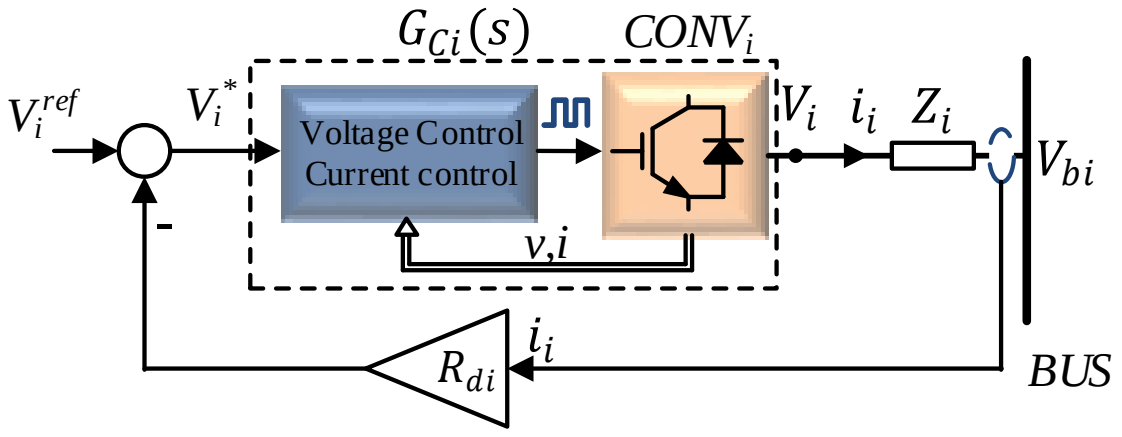

Fig. 2. Droop control implementation for converter $i$ ($CONV_i$).

*2.1. LTV Models for DC Microgrids with Conventional Droop Control*

As mentioned, the two main objectives, current sharing and bus voltage stability, are affected by the droop resistances (virtual resistances). Therefore, it is essential to review the droop control and formulate the relationships between (a) the virtual resistances and current sharing, and (b) virtual resistances and bus voltage in DC microgrids to support the adaptive droop design.

*2.1.1. Droop control*

Consider an example of droop control applied to a DC microgrid, which contains parallel DR (DC sources) $V_1$ and $V_2$ with initial reference voltages $V_1^{ref}$ and $V_2^{ref}$ connecting to the DC bus to supply power to a load $Z$ (Fig. 1). The terminal bus voltages are $V_{b1}$ and $V_{b2}$. Typically, the initial reference voltages are equal as $V_1^{ref} = V_2^{ref}$. Therefore, the steady-state relationship between the current sharing and the droop resistances is

$$\frac{i_1}{i_2} = \frac{R_{d2} + r_2}{R_{d1} + r_1}, \tag{1}$$

where $i_1$ and $i_2$ are the supplying currents, $r_1$ and $r_2$ are the equivalent line resistances from each DR to the load, and $R_{d1}$ and $R_{d2}$ are the droop resistances. The droop resistances are used to create a new setpoint for the output voltage of each DR to control the current sharing between resources.

*2.1.2. DC microgrid model*

Derivation of the DC microgrid model requires modeling of the three main components including DR, distribution systems, and loads. Take the system in Fig. 1 as an example. To consider the dynamics of DR, the droop control in Fig. 1 for one DC source is supposed to be carried out via one power converter ($CONV_i$) as shown in Fig. 2. The first component (DR) required for modeling is represented as

$$V_i = \left(V_i^{ref} - R_{di} i_i\right) G_{Ci}(s), \tag{2}$$

where the subscript ($i$) refers to the index of DR in the system ($i = 1, 2$), $V_i^{ref}$ is initial voltage references of DR$_i$, and $G_{Ci}(s)$ is the transfer function describing the relationships between the voltage reference and the output voltage of the DR. The second component, the interconnecting cables, is expressed as a series combination of resistance and inductance with the following relationship:

$$V_i - V_{bi} = r_i i_i + l_i \frac{di_i}{dt}. \tag{3}$$

The last component, a load with current consumption $i_L$, is related to the DR current $i_i$ as

$$i_L = \sum_{i=1}^{2} i_i \,. \tag{4}$$

To see the effect of droop changes on current sharing and bus voltage, a variation of the droop parameter $R_{di}$ by a small quantity $\hat{r}_{di}$ is conducted. Applying the small-signal analysis for (2) and (3), one obtains:

$$\begin{aligned} V_i + \hat{v}_i &= G_{Ci} V_i^{ref} - G_{Ci}(R_{di} + \hat{r}_{di})(i_i + \hat{\imath}_i) \\ V_i + \hat{v}_i - V_{bi} - \hat{v}_{bi} &= \left(r_i + \frac{d}{dt} l_i\right)(i_i + \hat{\imath}_i), \end{aligned} \tag{5}$$

where $\hat{v}_i$, $\hat{v}_{bi}$, and $\hat{\imath}_i$ are the small-signals of $V_i$, $V_{bi}$, and $i_i$, respectively. The relationships are simplified as

$$\begin{aligned} \hat{v}_i &= -G_{Ci}(\hat{r}_{di} i_i + R_{di}\hat{\imath}_i + \hat{r}_{di}\hat{\imath}_i) \\ \hat{v}_i &= \hat{v}_{bi} + \left(r_i + \frac{d}{dt} l_i\right) i_i. \end{aligned} \tag{6}$$

The parameters $\hat{\imath}_i$ and $\hat{r}_{di}$ are small values, which makes $\hat{r}_{di}\hat{\imath}_i \approx 0$. Additionally, two considerations are made for the small-signal analysis in (6) as follows:

(a) To formulate the small-signal relationship between bus voltage $\hat{v}_{bi}$ and the droop variation $\hat{r}_{di}$, the first consideration is made such that if the variation of current $\hat{\imath}_i$ is ignored ($\hat{\imath}_i = 0$) and $\hat{r}_{di}$ is taken as a main control signal, (6) becomes:

$$\begin{aligned} \hat{v}_i &= -G_{Ci}\hat{r}_{di} i_i \\ \hat{v}_i &= \hat{v}_{bi}. \end{aligned} \tag{7}$$

Thus, the following transfer function for the relationship between the bus voltage $\hat{v}_{bi}$ and the droop constant variation $\hat{r}_{di}$ is expressed as

$$\left.\frac{\hat{v}_{bi}}{\hat{r}_{di}}\right|_{\hat{\imath}_i = 0} = -i_i G_{Ci} \,. \tag{8}$$

(b) To formulate the small-signal relationship between current sharing $\hat{\imath}_i$ and the droop variation $\hat{r}_{di}$, the second consideration is made such that if in (6), the variation $\hat{v}_{bi}$ is ignored ($\hat{v}_{bi} = 0$) and $\hat{r}_{di}$ is taken as a main control signal one achieves:

$$\begin{aligned} \hat{v}_i &= -G_{Ci}(\hat{r}_{di} i_i + R_{di}\hat{\imath}_i) \\ \hat{v}_i &= \left(r_i + \frac{d}{dt} l_i\right)\hat{\imath}_i. \end{aligned} \tag{9}$$

Thus, the following transfer function for the relationship between the current sharing $\hat{\imath}_i$ and the droop constant $\hat{r}_{di}$ is expressed as

$$\left.\frac{\hat{\imath}_i}{\hat{r}_{di}}\right|_{\hat{v}_{bi} = 0} = -i_i G_{Ci} \frac{1}{R_{di} + r_i + s l_i}. \tag{10}$$

It is important to note that the small-signal relationships in (8) and (10) reflect the dynamic behavior of the bus voltage and current sharing on the droop resistance variation. These time-varying relationships lead to the adaptive reference models utilized in the next part for critical analysis of an advanced adaptive control system, and to a proposed adaptive droop secondary control system for DC microgrids.

*2.2. Development of Adaptive Control Algorithm*

To design an adaptive droop scheme for the DC microgrids control, we investigate a new trend of adaptive control called closed-loop reference model (CRM) control. The normalization technique paired with the projection algorithm is applied for CRM. Later, an adaptation gain-scheduling technique is proposed for various operating points of the control system.

#### 2.2.1. *Closed-loop reference model control*

Consider a general first-order system:

$$\dot{x} = ax + bu, \tag{11}$$

where $x$ and $u$ are the state variable and the control input; and $a$ and $b$ are unknown parameters, but the sign of $b$ is known. The control design criterion is to generate the bounded control signal $u$, by which the measured state $x$ follows the state $x_m$ of a stable reference model, which is defined as [38]:

$$\dot{x}_m = a_m x_m + b_m r - le, \tag{12}$$

where $r$ is the reference signal of the control system; $a_m$ and $b_m$ are known parameters with $a_m < 0$; $e = x - x_m$ is the model error; and $l < 0$ is the feedback parameter. The control update is selected as

$$u = \theta(t)x + k(t)r \tag{13}$$

where $\theta(t)$ and $k(t)$ are the controller parameters. Define the following vectors $\bar{\theta} = [\theta \quad k]^T$, $\phi = [x \quad r]^T$. The adaptation law is chosen as

$$\dot{\bar{\theta}} = -sgn(b)\gamma\phi e, \tag{14}$$

where $\gamma$ is the adaptation gain and $\gamma > 0$.

#### 2.2.2. *Proposed robust closed-loop reference model control*

The major focus of this part is on robust adaptive techniques including a normalization technique and a projection algorithm to improve the robustness of the CRM under model uncertainty and noise disturbances. Various setpoints are then studied to develop an advanced CRM controller using an adaptation gain-scheduling technique.

*a) Normalization technique*

The modeling error and signal vector $\phi$ are not guaranteed bounded. For example, the output of an unstable plant is unbounded. Hence, [33] introduces the normalization technique for the conventional adaptive method to guarantee that the normalized modeling error term and signal vector are bounded. Similarly, dynamic normalization is applied in this paper to the closed-loop adaptive method.
The error $e$ can be derived as

$$e = b^* \left(-\bar{\theta}^{*T}\phi_n + u_n\right), \tag{15}$$

where $b^*$ is the expected value of $b$, $\bar{\theta}^*$ is the expected value of $\bar{\theta}$, $\phi_n = \frac{1}{s-a_m-l}\phi$, and $u_n = \frac{1}{s-a_m-l}u$. The estimation error $\hat{e}$ of $e$ is expressed as

$$\hat{e} = b(-\bar{\theta}^T\phi_n + u_n), \tag{16}$$

where $b$ is the online estimate of $b^*$. Then, the modeling error is defined as

$$\varepsilon = \frac{e - \hat{e}}{m^2}, \tag{17}$$

where $m = \sqrt{1 + {\phi_n}^T\phi_n + {u_n}^2}$ is the normalized term to create the boundedness of $\varepsilon$, $\frac{\phi_n}{m}$, and $\frac{u_n}{m}$. The adaptation law is chosen as

$$\begin{aligned}\dot{\bar{\theta}} &= -\gamma sgn(b)\phi_n\varepsilon \\ \dot{b} &= \gamma(-\bar{\theta}^T\phi_n + u_n)\varepsilon.\end{aligned} \tag{18}$$

Therefore, the above parameters $\dot{\bar{\theta}}$ and $\dot{b}$ are bounded because the quantities $\frac{\phi_n}{m}, \frac{u_n}{m}, \frac{e}{m}, \frac{\hat{e}}{m}$, and $\varepsilon$ are bounded.

*b) Parameter projection*

In adaptive control, bounded disturbances and errors in dynamic models can cause the control parameters to drift to infinity, which results in system's instability. Therefore, a common modification called parameter projection in the adaptive law is introduced to prevent the parameter drift phenomenon. The modification is made based on some known properties of the system. Suppose that the knowledge about the system results in the statement that the controller parameters ($\bar{\theta}$ and $b$) could only drift to the boundary of the following known convex sets:

$$\begin{aligned}M_{\bar{\theta}} &= \{\bar{\theta} | g(\bar{\theta}) = \bar{\theta}^T\bar{\theta} - M_{0\bar{\theta}}^2 \le 0\} \\ M_b &= \{b | g(b) = b^2 - M_{0b}^2 \le 0\},\end{aligned} \tag{19}$$

where $M_{0\bar{\theta}}$ and $M_{0b}$ are known parameters; and $g(\bar{\theta})$ and $g(b)$ are functions to describe the constraint of control parameters. Thus, the following parameter projection algorithm is performed for parameter drift prevention as

$$\dot{\bar{\theta}} = \begin{cases} \dot{\bar{\theta}}_0 & \text{if } \bar{\theta}^T\bar{\theta} \le M_{0\bar{\theta}}^2 \text{ or} \\ & \text{if } \bar{\theta}^T\bar{\theta} = M_{0\bar{\theta}}^2 \text{ and } (\dot{\bar{\theta}}_0)^T\bar{\theta} \le 0 \\ \left(I - \dfrac{\gamma\bar{\theta}\bar{\theta}^T}{\bar{\theta}^T\gamma\bar{\theta}}\right)\dot{\bar{\theta}}_0 & \text{Otherwise} \end{cases}$$

$$\dot{b} = \begin{cases} \dot{b}_0 & \text{if } b^2 \le M_{0b}^2 \text{ or} \\ & \text{if } b^2 = M_{0b}^2 \text{ and } \dot{b}_0 b \le 0 \\ 0 & \text{Otherwise,} \end{cases} \tag{20}$$

where $\dot{\bar{\theta}}_0 = -\gamma sgn(b)\phi_n\varepsilon$, $\dot{b}_0 = \gamma(-\bar{\theta}^T\phi_n + u_n)\varepsilon$, and $I_{2\times2}$ is the $(2 \times 2)$ identity matrix.

*c) Adaptation gain-scheduling*

Since, there is a significant impact of the reference points on the system response, an adaptation gain-scheduling technique is proposed. To observe the impact of the adaptive gain on the system's response due to the change in the operating point, the control update is considered for analysis as

$$u = \phi^T\bar{\theta} = \phi^T \int -\gamma sgn(b)\phi_n\varepsilon \, dt. \tag{21}$$

In the discrete-time domain,

$$\Delta u_k = u_{k+1} - u_k = \phi_k^T\phi_{nk}[-\gamma sgn(b)\varepsilon_k T_s], \tag{22}$$

where $T_s$ is the sampling time and $k$ is the time instant index. (22) shows that if there is a change in the reference signal $r$ by a ratio $\alpha$, the new reference $\alpha r$ will approximately change the state response to $\alpha x$. These changes lead to the following change in the control input:

$$\Delta u_k = u_{k+1} - u_k = \alpha^2\phi_k^T\phi_{nk}[-\gamma sgn(b)\varepsilon_k T_s]. \tag{23}$$

Thus, the term $\phi = [x \quad r]^T$ has a major impact on the control input $u$ because $\Delta u_k$ changes by the ratio of $\alpha^2$ when the system setpoint $r$ changes. To maintain the control performance in various operating conditions, the following modifications are made for the incremental control input $\Delta u_k$ and the adaptation gain $\gamma_k$:

$$\begin{aligned}\Delta u_k &= u_{k+1} - u_k = \gamma_k\phi_k^T\phi_{nk}[-sgn(b)\varepsilon_k T_s] \\ \gamma_k &= \gamma_0\frac{1}{\alpha_k^2},\end{aligned} \tag{24}$$

where $\gamma_k$ is the adaptation gain at time instant $k$, $\gamma_0$ is the initial design adaptation gain for the initial desired reference signal ($r_0 = 1$), and $\alpha_k = \frac{r_k}{r_0} = r_k$. Therefore, the adaptation law is selected as

$$\dot{\bar{\theta}} = \begin{cases} \dot{\bar{\theta}}_0 & \text{if } \bar{\theta}^T\bar{\theta} \le M_{0\bar{\theta}}^2 \text{ or} \\ & \text{if } \bar{\theta}^T\bar{\theta} = M_{0\bar{\theta}}^2 \text{ and } (\dot{\bar{\theta}}_0)^T\bar{\theta} \le 0 \\ \left(I - \dfrac{\gamma_k\bar{\theta}\bar{\theta}^T}{\bar{\theta}^T\gamma_k\bar{\theta}}\right)\dot{\bar{\theta}}_0 & \text{Otherwise} \end{cases}$$

$$\dot{b} = \begin{cases} \dot{b}_0 & \text{if } b^2 \le M_{0b}^2 \text{ or} \\ & \text{if } b^2 = M_{0b}^2 \text{ and } \dot{b}_0 b \le 0 \\ 0 & \text{Otherwise,} \end{cases} \tag{25}$$

where $\dot{\bar{\theta}}_0 = -\gamma_k sgn(b)\phi_n\varepsilon$ and $\dot{b}_0 = \gamma_k(-\bar{\theta}^T\phi_n + u_n)\varepsilon$.

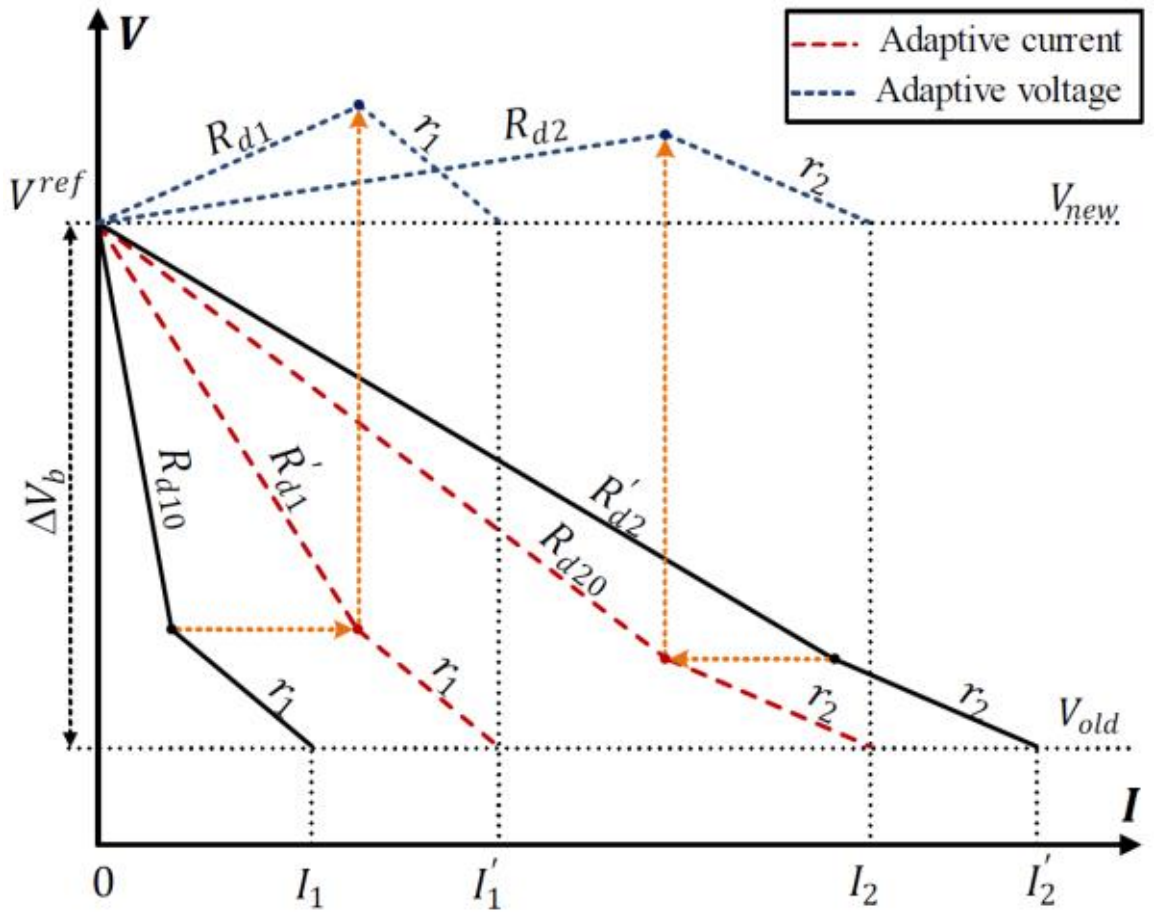

Fig. 3. Two-stage adaptive droop characteristics.

**Theorem 1.** *The adaptation law (25) with control update (13) guarantees that the system (11) is stable and the system state $x$ tracks the reference model state $x_m$ asymptotically.*

The proof of **Theorem 1** is illustrated in the appendix.

### *2.3. Proposed Adaptive Droop Control*

#### *2.3.1. Proposed adaptive droop control architecture*

The idea of the proposed adaptive droop is illustrated through the droop diagram in Fig. 3. The initial droop characteristic lines $R_{d10} + r_1$ and $R_{d20} + r_2$ result in non-desired current sharing $I_1$ and $I_2$, and bus voltage deviation $\Delta V_b$. The quantities $I'_1$ and $I'_2$ represent the desired current sharing and $V_{new}$ represents the desired voltage operation level. The adaptive control is applied to change the droop characteristic lines to a desired position, which satisfies the DC microgrid control criteria. The droop changes can be illustrated in two stages. The first stage is the adaptive current control, where the initial droop lines are moved horizontally to the new positions $R'_{d1} + r_1$ and $R'_{d2} + r_2$ for the desired current sharing $I'_1$ and $I'_2$. The second stage minimizes the bus voltage deviation by moving the droop lines vertically to the new positions $R_{d1} + r_1$ and $R_{d2} + r_2$, where $V_{new} = V^{ref}$. As a result, the final droop value $R_{di}$ for $\text{DR}_\text{i}$ $(i = 1, 2, \dots, N)$ that satisfies the current sharing and bus voltage control in DC microgrids is calculated as follows:

$$R_{di} = R_{di0} + R_{Vi} + R_{Ii}, \tag{26}$$

where $R_{di0}$ is the initial droop value, $R_{Vi}$ is the adaptive droop component for bus voltage restoration, and $R_{Ii}$ is the adaptive droop component for current sharing control.

The adaptive droop resistances $R_{Vi}$ and $R_{Ii}$ are implemented in a distributed secondary control block $\text{DSC}_\text{i}$ of node $i$ as shown in Fig. 4. There are three layers in the distributed architecture. The first layer is the primary control, where droop control and the internal voltage and current control of a power converter are performed. The second layer is the secondary control, where adaptive control is conducted for current sharing and bus voltage stability. The third layer is the communication layer for exchanging the control information including bus voltage and current sharing information. In detail, the $\text{DSC}_i$ sends its terminal per-unit voltage $v_{i,pu}$ and per-unit current $i_{i,pu}$ information and receives per-unit voltage $v_{j,pu}$ and per-unit current $i_{j,pu}$ information from its neighbor $\text{DSC}_j$ via the communication layer.

In this architecture, the communication network is represented as a directed graph $G = (\mathcal{V},\ \mathcal{E})$. $\mathcal{V}$ denotes the set of nodes, where $\mathcal{V} =$ $(1, 2, \dots, i, \dots, N)$. $\mathcal{E}$ represents the set of edges, where $\mathcal{E} \subset \mathcal{V} \times \mathcal{V}$. Neighbors of $i$ belong to a set defined as $N_i = \{j \in \mathcal{V}\ |j \neq i, (j, i) \in \mathcal{E}\}$. The adjacency matrix $A = [a_{ij}]$ is defined such that $a_{ij} = 1$ if $(j, i) \in \mathcal{E}$, and $a_{ij} = 0$ otherwise. Based on that, the reference current $i_{i,pu}^{ref}$ and the average voltage $\bar{v}_{i,pu}$ are calculated as

$$i_{i,pu}^{ref} = \frac{\sum_{j=1}^{N_j} a_{ij} i_{j,pu}}{\sum_{j=1}^{N_j} a_{ij}}$$
$$\bar{v}_{i,pu} = \frac{v_{i,pu} + \sum_{j=1}^{N_j} a_{ij} v_{j,pu}}{1 + \sum_{j=1}^{N_j} a_{ij}}, \tag{27}$$

which are used as the setpoint for current sharing control and feedback value for bus voltage restoration.

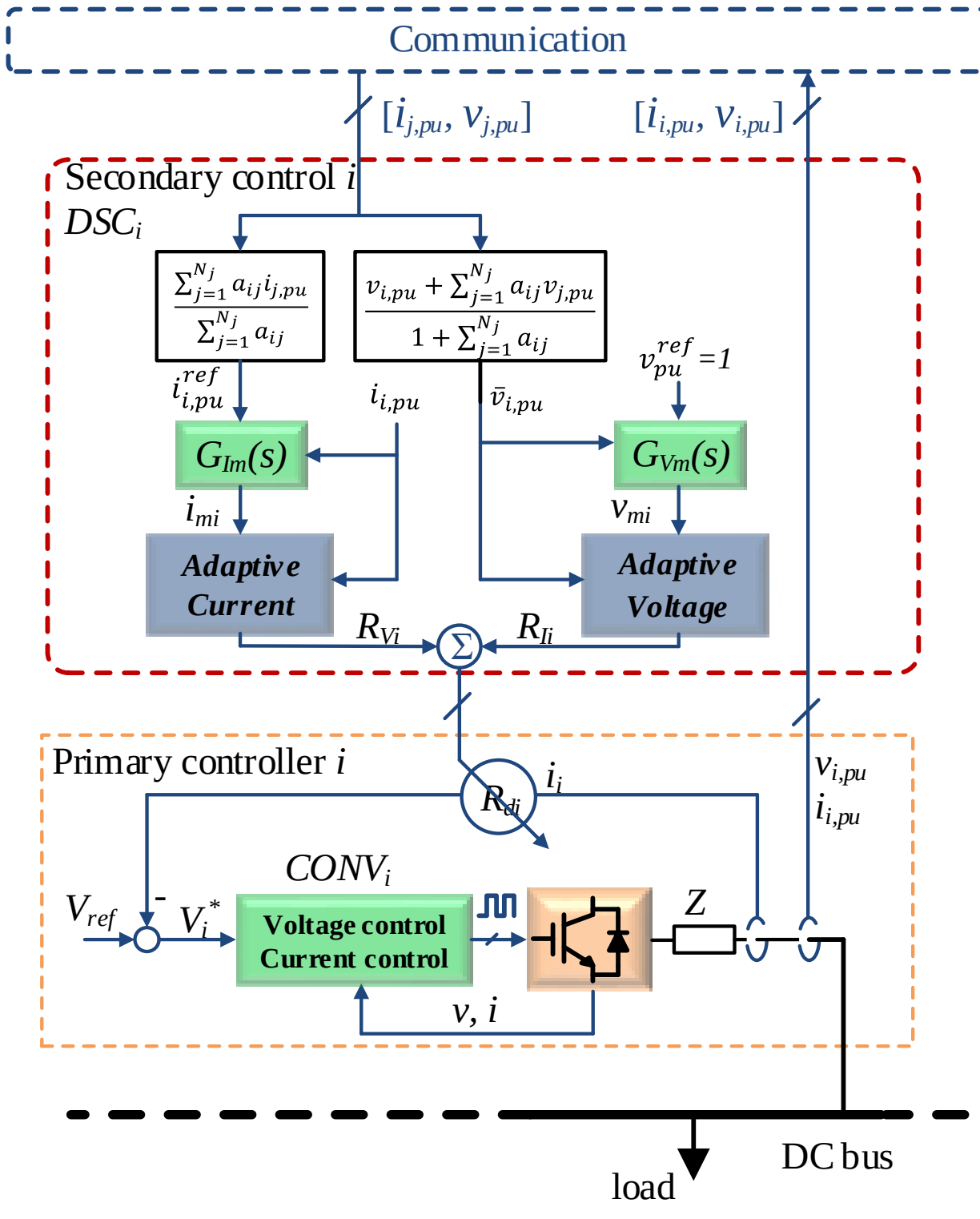

Fig. 4. Adaptive droop secondary control for DC microgrids.

Regarding control implementation, the reference current $i_{i,pu}^{ref}$ and feedback current $i_{i,pu}$ are the inputs for the adaptive current controller. The selected reference bus voltage $v_{pu}^{ref} = 1$ and average bus voltage $\bar{v}_{i,pu}$ are the inputs of the distributed adaptive voltage controller.

#### *2.3.2. Voltage and current adaptation mechanism design*

The relationships between the current or voltage and the droop resistances (8), (10) represent the system models. In these models, the fast dynamics of power converters with multiple-order transfer functions can be neglected. Thus, the voltage transfer function $G_{Ci}$ for the converter $i$ is assumed as the following first-order system:

$$G_{Ci} = \frac{1}{1 + \tau_{vi} s}. \tag{28}$$

Therefore (8), (10) become

$$G_{Vi}(s) = \frac{\hat{v}_{bi}}{\hat{r}_{di}} = -i_i \frac{1}{1 + \tau_{vi} s}$$
$$G_{Ii}(s) = \frac{\hat{\imath}_i}{\hat{r}_{di}} \cong -\frac{i_i}{R_{di} + r_i} \frac{1}{1 + \left(\frac{l_i}{R_{di} + r_i} + \tau_{vi}\right) s}. \tag{29}$$

By relating the small-signals (a) $\hat{v}_{bi}$ to the real average bus voltage $\bar{v}_{i,pu}$, (b) $\hat{\imath}_i$ to the real current $i_{i,pu}$, and (c) $\hat{r}_{di}$ to the resistance quantities $R_{Vi}$ and $R_{Ii}$, the following models are achieved:

$$G_{Vi}(s) = \frac{\bar{v}_{i,pu}}{R_{Vi}} = -i_i \frac{1}{1+\tau_{vi}s}$$
$$G_{Ii}(s) = \frac{i_{i,pu}}{R_{Ii}} \cong -\frac{i_i}{R_{di}+r_i} \frac{1}{1+\left(\frac{l_i}{R_{di}+r_i}+\tau_{vi}\right)s}, \tag{30}$$

which can be rewritten as

$$G_{Vi}(s) = \frac{\bar{v}_{i,pu}}{R_{Vi}} = \frac{b_{Vi}}{s-a_{Vi}}$$
$$G_{Ii}(s) = \frac{i_{i,pu}}{R_{Ii}} = \frac{b_{Ii}}{s-a_{Ii}}, \tag{31}$$

where $b_{Ii}$, $a_{Ii}$, $b_{Vi}$, and $a_{Vi}$ are unknown, but the signs of these parameters are known as $sgn(b_{Ii}) = -sgn(i_i)$, $sgn(b_{Vi}) = -sgn(i_i)$, $a_{Ii} < 0$, $a_{Vi} < 0$. The relationships in (31) represent the reduced first-order models of the system. As mentioned in the Introduction, it is not necessary to investigate the full-order model of the system since the adaptive control process does not require the in-depth knowledge of the system. The control parameters will be identified online based on the known information via the system models and based on specified desired reference characteristics the system should follow.

The unknown first-order systems (31) lead to a decision that $G_{Vmi}(s)$ and $G_{Imi}(s)$ shown in (32) are selected as the bus voltage and current reference model for the proposed adaptive control design of $DSC_i$.

$$G_{Vmi}(s) = \frac{v_{mi}}{v_{pu}^{ref}} = \frac{b_{Vmi}}{s-a_{Vmi}}$$
$$G_{Imi}(s) = \frac{i_{mi}}{i_{i,pu}^{ref}} = \frac{b_{Imi}}{s-a_{Imi}}, \tag{32}$$

where $v_{pu}^{ref}$ and $i_{i,pu}^{ref}$ are the setpoints; $v_{mi}$ and $i_{mi}$ are the bus voltage and current output of the reference models; and $b_{Vmi} > 0$, $a_{Vmi} < 0$, $b_{Imi} > 0$, and $a_{Imi} < 0$ are the desired parameters of the reference model.

**Theorem 1** is utilized for the adaptive bus voltage and current sharing control. To set up the adaptive algorithm for bus voltage and current sharing control, the following normalized terms for converter $i$ are defined:

$$\phi_{Vni} = \frac{1}{s-a_{Vmi}-l_{Vi}}\left[\bar{v}_{i,pu} \quad v_{pu}^{ref}\right]^T$$
$$R_{Vni} = \frac{1}{s-a_{Vmi}-l_{Vi}} R_{Vi}$$
$$\phi_{Ini} = \frac{1}{s-a_{Imi}-l_{Ii}}\left[i_{i,pu}^{ref} \quad i_{i,pu}\right]^T$$
$$R_{Ini} = \frac{1}{s-a_{Imi}-l_{Ii}} R_{Ii}. \tag{33}$$

Thus, the adaptation mechanism for bus voltage control is formulated as

$$R_{Vi} = \bar{\theta}_{Vi}^T\left[\bar{v}_{i,pu} \quad v_{pu}^{ref}\right]^T$$
$$\dot{\bar{\theta}}_{Vi} = -\gamma_{Vik}\, sgn(b_{Vi})\phi_{Vni}\varepsilon_{Vi}$$
$$\dot{b}_{Vi} = \gamma_{Vik}(-\bar{\theta}_{Vi}^T\phi_{Vni} + R_{Vni})\varepsilon_{Vi}$$
$$\varepsilon_{Vi} = \frac{e_{Vi}-\hat{e}_{Vi}}{m_{Vi}^2}, \tag{34}$$

where $e_{Vi} = \bar{v}_{i,pu} - v_{mi}$, $\hat{e}_{Vi} = b_{Vi}(-\bar{\theta}_{Vi}^T\phi_{Vni} + R_{Vni})$, $m_{Vi} = \sqrt{1+\phi_{Vni}^T\phi_{Vni}+R_{Vni}^2}$, $\gamma_{Vik} = \frac{\gamma_{Vi0}}{\alpha_{Vik}^2}$, and $\alpha_{Vik} = v_{pu}^{ref}$. Similar to the voltage control, the current adaptive control mechanism is designed as

$$R_{Ii} = \bar{\theta}_{Ii}^T\left[i_{i,pu} \quad i_{i,pu}^{ref}\right]^T$$
$$\dot{\bar{\theta}}_{Ii} = -\gamma_{Iik}\, sgn(b_{Ii})\phi_{Ini}\varepsilon_{Ii}$$
$$\dot{b}_{Ii} = \gamma_{Iik}(-\bar{\theta}_{Ii}^T\phi_{Ini} + R_{Ini})\varepsilon_{Ii}$$
$$\varepsilon_{Ii} = \frac{e_{Ii}-\hat{e}_{Ii}}{m_{Ii}^2}, \tag{35}$$

where $e_{Ii} = i_{i,pu} - i_{mi}$, $\hat{e}_{Ii} = b_{Ii}(-\bar{\theta}_{Ii}^T\phi_{Ini} + R_{Ini})$, $m_{Ii} = \sqrt{1+\phi_{Ini}^T\phi_{Ini}+R_{Ini}^2}$, $\gamma_{Iik} = \frac{\gamma_{Ii0}}{\alpha_{Iik}^2}$, and $\alpha_{Vik} = i_{i,pu}^{ref}$.

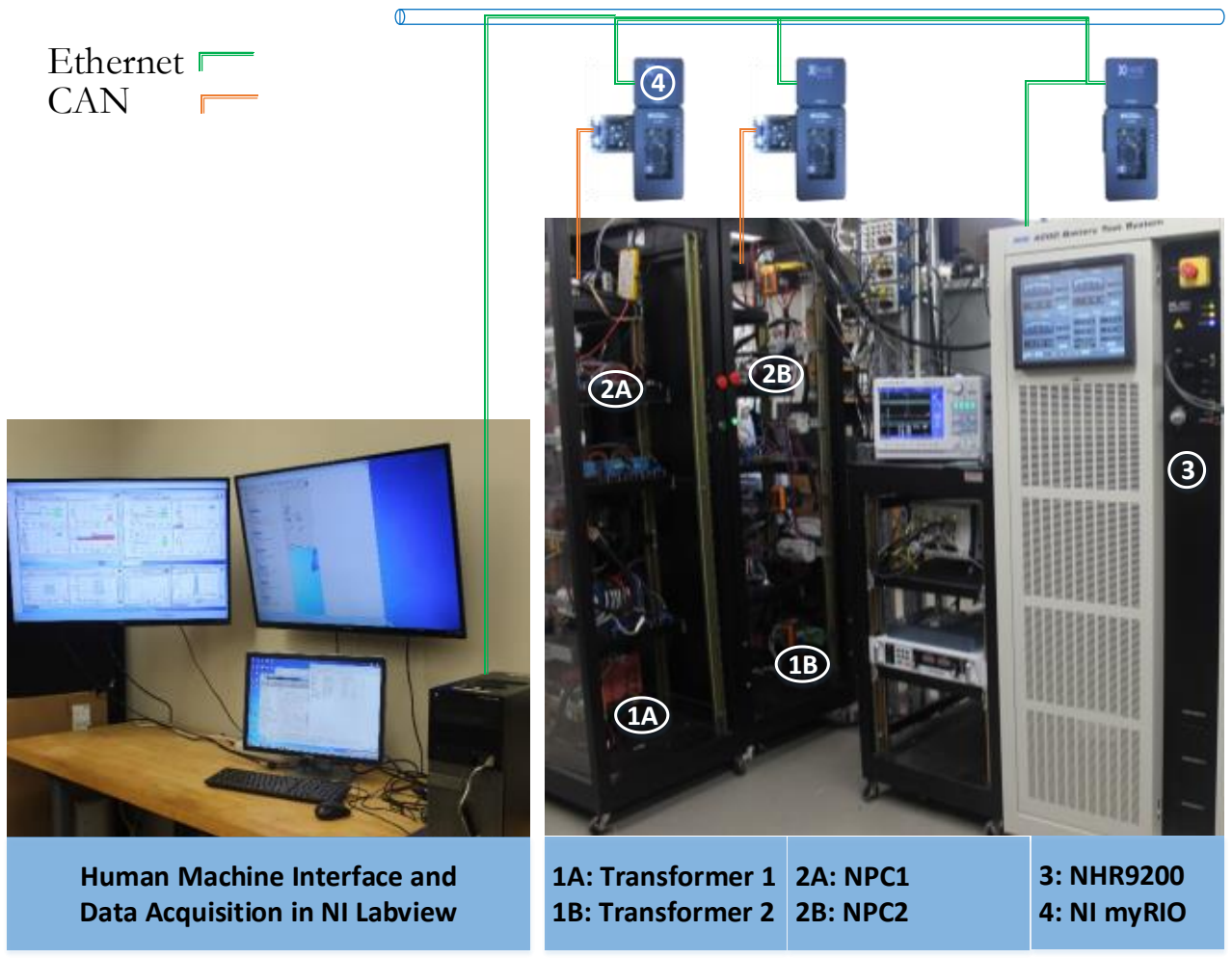

Fig. 5. Experimental system.

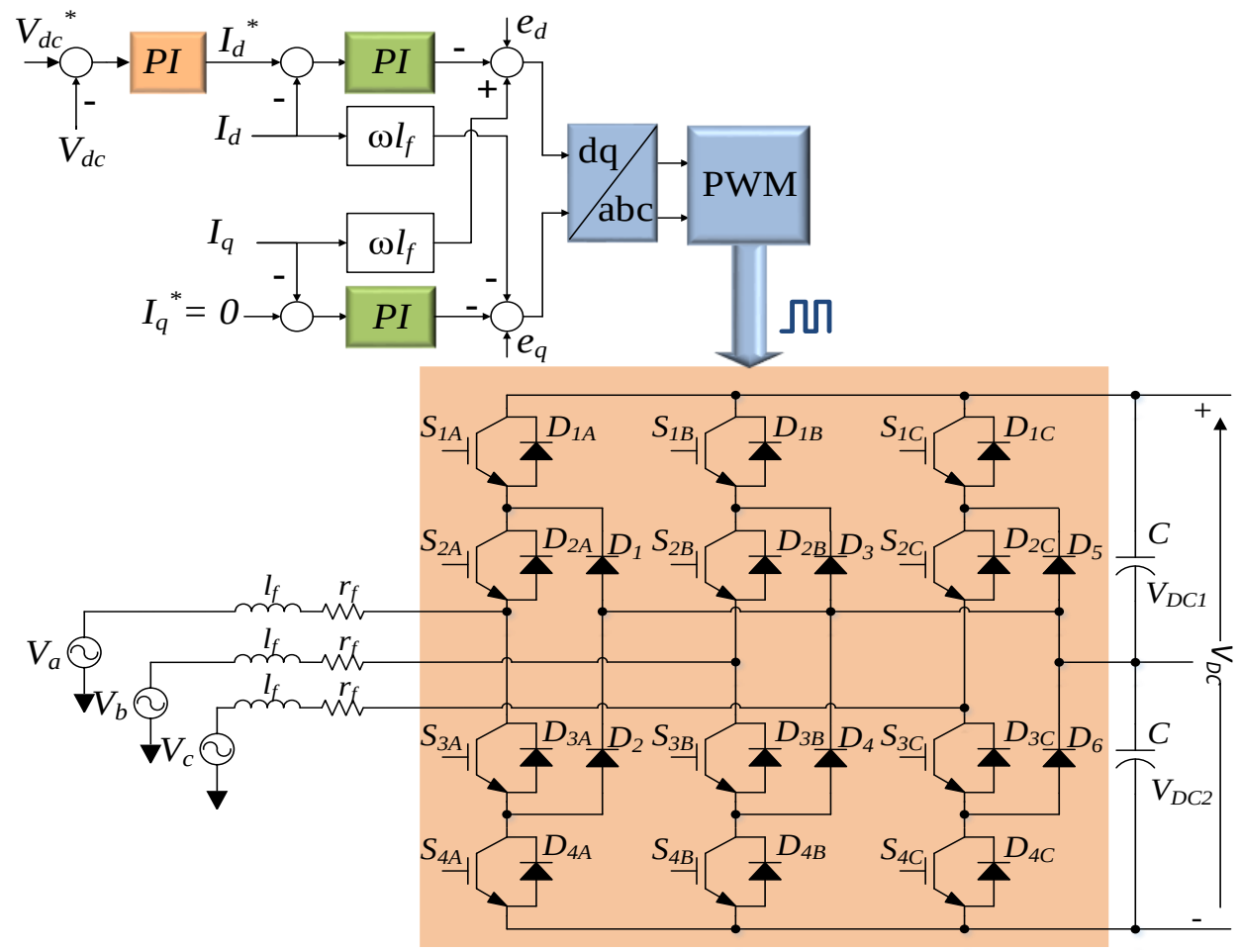

Fig. 6. NPC topology and control diagram.

Additionally, to prevent control parameter ($\bar{\theta}_{Vi}, b_{Vi}, \bar{\theta}_{Ii}$, and $b_{Ii}$) drift when updating them using (34) and (35), the projection algorithm (24) is applied. Based on that, the droop resistances $R_{Vi}$ and $R_{Ii}$ are identified iteratively. Consequently, the final droop for a converter $i$ is updated according to (26).

## 3. Case Studies and Results

### 3.1. System Description

The experimental setup is shown in Fig. 5. The microgrid used to verify the control algorithm includes two three-level neutral point clamped converters (NPC) operated as rectifiers. The two NPC are powered by two AC transformers, which have the same AC source from the laboratory grid. The NPC converters are controlled in the $dq$ reference frame as shown in Fig. 6. The DC output voltage $V_{dc}$ of the NPC is controlled with a PI controller to generate the current reference $I_d^*$. The other reference current $I_q^*$, is set to 0 to achieve the

TABLE I

DC Microgrid Parameters

| Symbol | Quantity | Values |
|---|---|---|
| $V_{L-L}$ | NPC input voltage | 208 V (60Hz) |
| $V_{DC}$ | NPC output voltage | 400 V |
| $l_i$ | NPC input filter inductor | 2.07 mH |
| C | NPC output capacitor | 380 µF |
| $r_{fi}$ | NPC input resistor | 0.2 Ω |
| $f$ | Switching frequency | 20 kHz |
| $R_i$ | Cable resistance | 0.5 Ω |
| $L_i$ | Cable inductance | 3 mH |
| $\tau_{vi}$ | NPC voltage loop time constants | 5 ms |
| $P_{NPC_1}$ | $NPC_1$ nominal power | 4 kW |
| $P_{NPC_2}$ | $NPC_2$ nominal power | 2 kW |

TABLE II

Adaptive Control Parameters

| Symbol | Quantity | Values |
|---|---|---|
| $a_{mVi}$ | Voltage controller parameter | -10 |
| $b_{mVi}$ | Voltage controller parameter | 10 |
| $a_{mIi}$ | Current controller parameter | -10 |
| $b_{mIi}$ | Current controller parameter | 10 |
| $\gamma_{Vi0}$ | Voltage adaptation gain | 1000 |
| $\gamma_{Ii0}$ | Current adaptation gain | 1000 |
| $l_{Vi}$ | Voltage feedback gain | −10 |
| $l_{Ii}$ | Current feedback gain | −10 |

unity power factor. Decoupling and compensation control methods using input inductor ($l_f$) with $dq$-input voltages ($e_d$ and $e_q$) are also utilized. More details on this converter control can be found in [44]. The two NPC provide 400 $V_{dc}$ to the DC bus, which supplies power to a load device (NHR9200). 4 kW and 2 kW are the rated power assumed for $NPC_1$ and $NPC_2$, respectively. Information about the hardware configuration is indicated in TABLE I. The NPC and load are controlled by the secondary controllers (NI myRIO), located above the DSP of the two NPC and the NHR9200.

The communication between the primary controller (DSP) of the NPC and secondary controller (NI myRIO) is conducted via CAN (1Mb/s), where the DSP sends current and voltage information to the NI myRIO, and the NI myRIO sends droop information back to the DSP. The communication among NI myRIO is implemented via Ethernet TCP/IP (100Mb/s), where the NI myRIO controllers exchange current and voltage information with each other based on an undirected communication graph. The TCP/IP communication among NI myRIO showed communication delays up to 10 ms; thus, this value was chosen as the time step for the adaptive system.

Because the dynamics of the two NPC are identical, the free design parameters of the adaptive voltage and adaptive current reference models for the two NPC are selected to be equal (TABLE II). The adaptive control parameters were selected based on the simulation analysis of the first-order system in (31) with reference model (32).

There are two test cases made to verify the effectiveness of the proposed methodology. The first case is the adaptive control activation, which is made while running the conventional droop method at 3 kW to prove that the adaptive method improves the limitations when compared to the conventional one. In this test case, the current sharing control and bus voltage restoration are activated sequentially. The second test case is made by changing the load power from 3 kW to 6 kW to ensure the stability of the method under load disturbance.

### *3.2. Simulation Results*

The simulations were conducted in MATLAB/Simulink. Results for current control activation, voltage control activation, and constant power load change from 3 kW to 6 kW are shown in Fig. 7, Fig. 8, and Fig. 9, respectively.

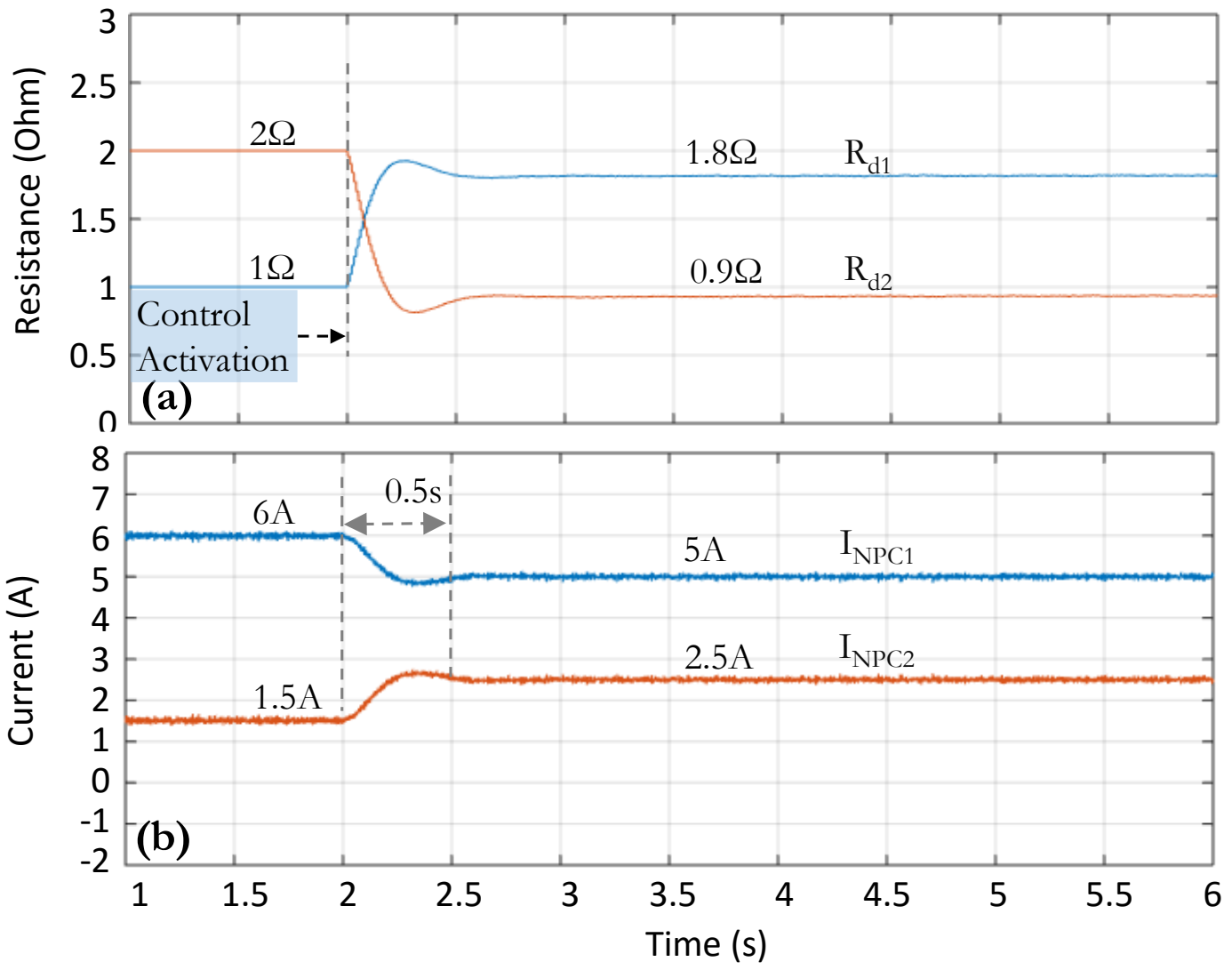

Fig. 7. Current sharing performance under control activation. (a) droop resistances, (b) output currents.

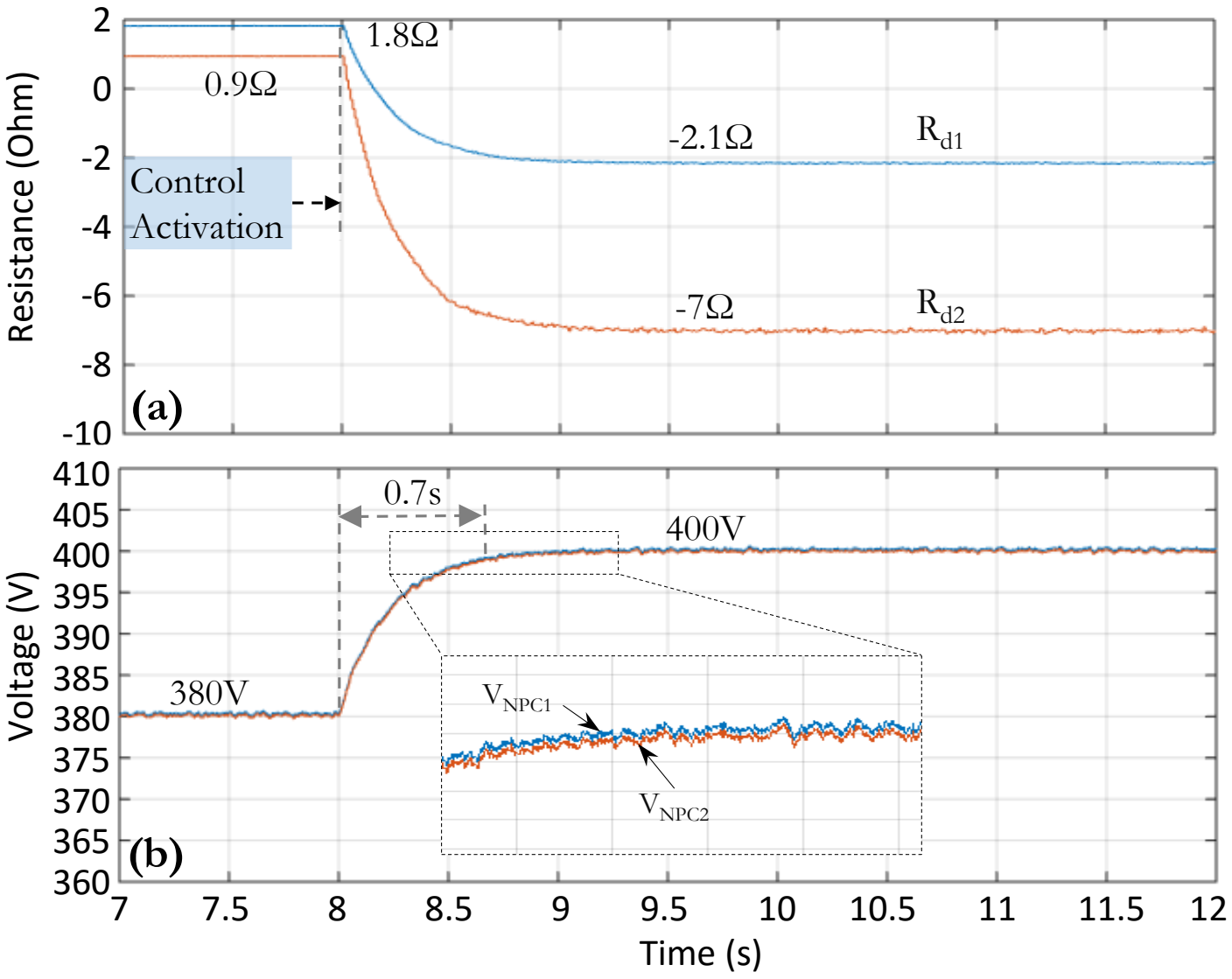

Fig. 8. Bus voltage restoration under control activation. (a) droop resistances, (b) bus voltages.

### *3.3. Experimental Results*

The proposed adaptive control algorithm was implemented in the NI myRIO controllers with a 10 ms sampling time step. Data acquisition in the LabVIEW graphical interface shows the adaptive droop profile during each stage of operation of the system (Fig. 10). The experimental results for current sharing and bus voltage restoration control activation are illustrated in Fig. 11a and Fig. 11b, respectively. Fig. 12 shows the transient in current sharing between $NPC_1$ and $NPC_2$, and the bus voltage restoration during an increase in load (NHR9200) power.

## 4. Analysis and Discussions

The requirements for the control system are to share the current between DR as 2:1 for $NPC_1$ and $NPC_2$ and to maintain the bus voltages at 400 V in steady state operation.

Simulation results for the first test case are shown in Fig. 7 and Fig. 8. The first stage of control activation is shown in Fig. 7, where the current among the NPC's are shared proportionally. The second stage of the control activation is shown in Fig. 8, where the bus voltage is recovered to the nominal value. In this test case, the load operates as a resistive load at 0.5 pu (7.5 A). The first stage results in the change of the droop resistances (Fig. 7a), which corresponds to the accurate current sharing (Fig. 7b). The conventional method operating before 2s with droop resistances $R_{d1} = 1\ \Omega$ and $R_{d2} = 2\ \Omega$ results in the non-proportional current sharing, which is 6 A for $NPC_1$ and 1.5 A for $NPC_2$. The adaptive current controller activated at 2s changes the droop resistances to 1.8 Ω and 0.9 Ω for $R_{d1}$ and $R_{d2}$, respectively. This activation process takes 0.5s to share the current proportionally with 5 A and 2.5 A for $NPC_1$ and $NPC_2$, respectively. In the second stage, the bus voltage control is activated at 8s. The adaptive voltage controller brings the droop resistances down to -2.1 Ω and -7 Ω for $R_{d1}$ and $R_{d2}$ (Fig. 8a). According to (2), these negative droops will increase the converters' terminal voltages to recover the DC bus voltage to the nominal value at 400 V. The droop resistances change takes 0.7s to restore the bus voltage from 380 V to the nominal value of 400 V (Fig. 8b). The experimental results in Fig. 10 for droop changes in current control activation are 1.8 Ω for $R_{d1}$ and 0.8 Ω for $R_{d2}$. For the bus voltage control activation, the droop changes are -2.1 Ω for $R_{d1}$ and -5 Ω for $R_{d2}$. These droop changes are slightly different from the droop changes in simulation because of the uncertainty in the experimental system. However, these changes in droop resistances share the current accurately and recover the bus voltages to 400 V. Fig. 11a indicates that the current sharing activation takes 0.5s to achieve the accurate current sharing, in which $NPC_1$ supplies 5 A and $NPC_2$ supplies 2.5 A to the NHR9200 load. Fig. 11b shows that the bus voltages are recovered from 380 V to 400 V in 0.7s after activating the adaptive voltage controller.

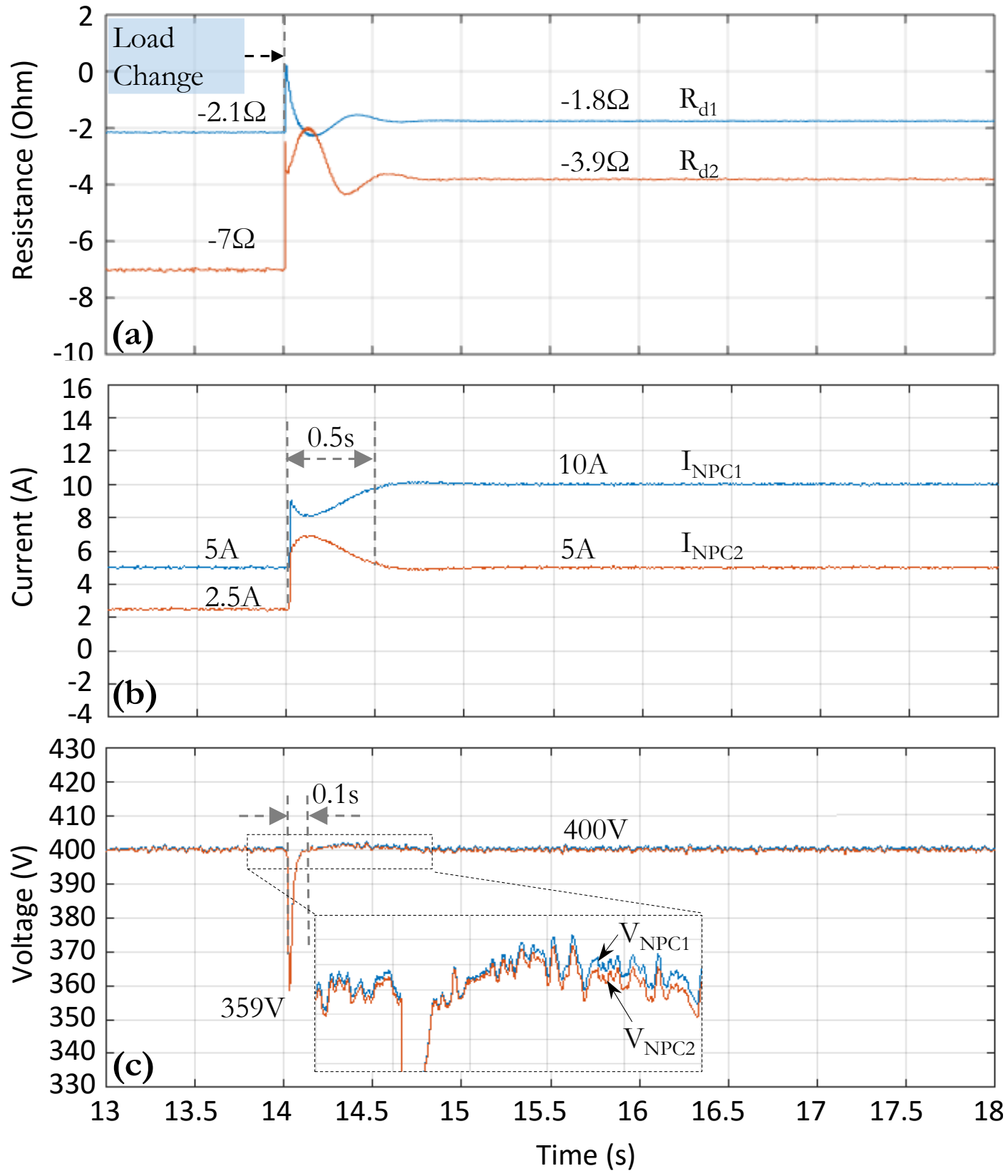

Fig. 9. Adaptive performance under load change. (a) droop resistances, (b) output currents, (c) bus voltages.

In the second test case, the nonlinear load changes from 0.5 pu (7.5 A) to 1 pu (15 A) at 14s. Simulation results of this case are shown in Fig. 9 with droop resistances presented in Fig. 9a, supply current presented in Fig. 9b, and the bus voltage regulation shown in Fig. 9c. Fig. 9a indicates that once the load changes, the droop resistances change to -1.8 Ω and -3.9 Ω for $R_{d1}$ and $R_{d2}$. The change in the droop resistances helps to maintain the proportional current sharing as $I_{NPC1} = 5\ A$ and $I_{NPC2} = 10\ A$ (Fig. 9b). The change in the droop resistances also helps to restore the bus voltage control to 400 V after the 0.1s transient time (Fig. 9c). In the experiment, the DC load (NHR9200) is set up to draw 3 kW more power to have a total power consumption at 6 kW. Results in Fig. 10 indicate that the droop operation points change to $R_{d1} = -1.1\ \Omega$ and $R_{d2} = -2.2\ \Omega$ to maintain the system's objectives. Fig. 12a illustrates that the adaptive control is able to maintain the desired proportional current sharing between two converters after 0.5s of transient time, where $NPC_1$ supplies 10 A and $NPC_2$ supplies 5 A to the NHR9200. In the transient time of this scenario, the bus voltage drops to 350 V, but is restored to 400 V in 0.1s by the adaptive control (Fig. 12b).

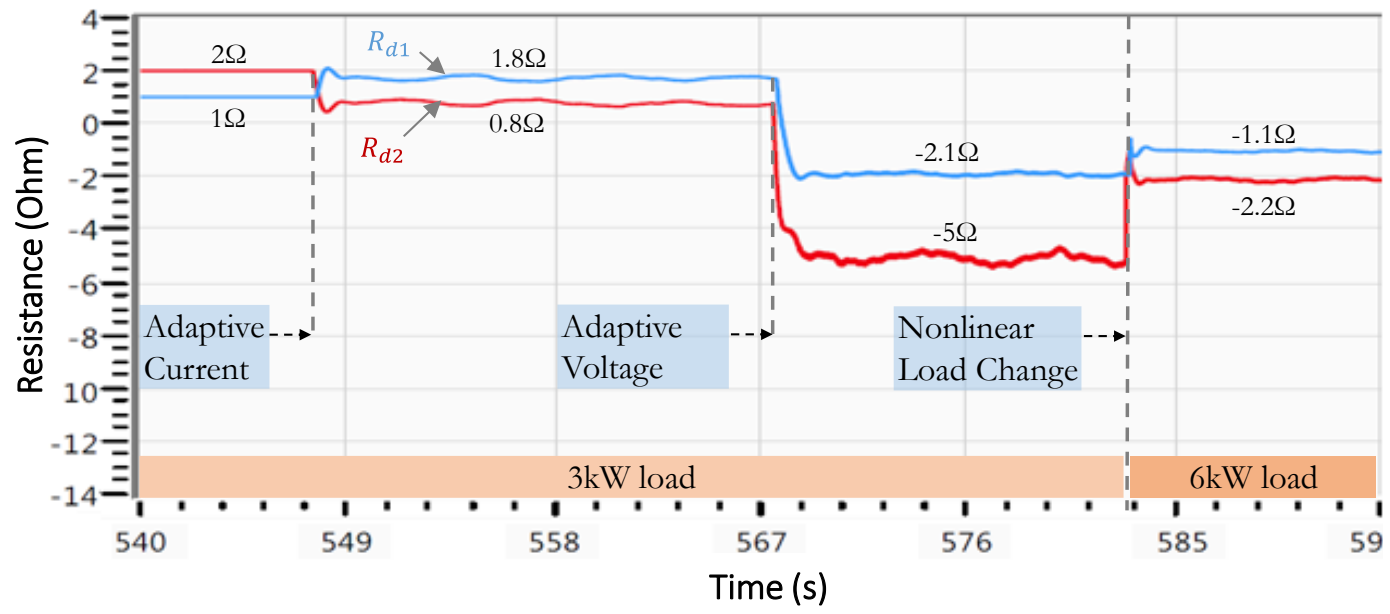

Fig. 10. Adaptive droop profile acquired in LabVIEW.

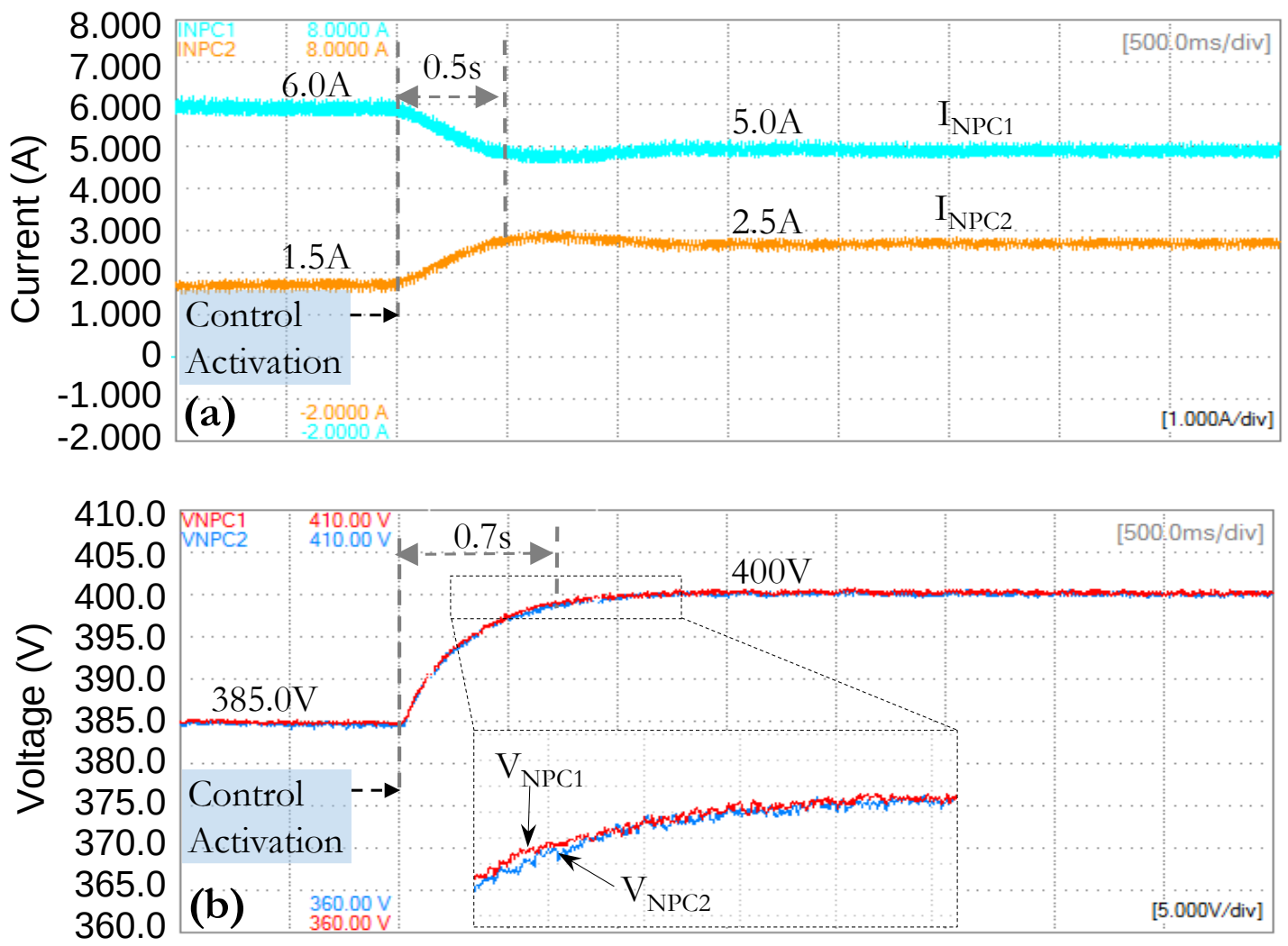

Fig. 11. Control Activation. (a) adaptive current, (b) adaptive voltage control. Voltage: 10V/div. Current: 1A/div. Horizontal axis: 500ms/div. Zoomed voltage area (Voltage: 1V/div, Horizontal axis: 200ms/div).

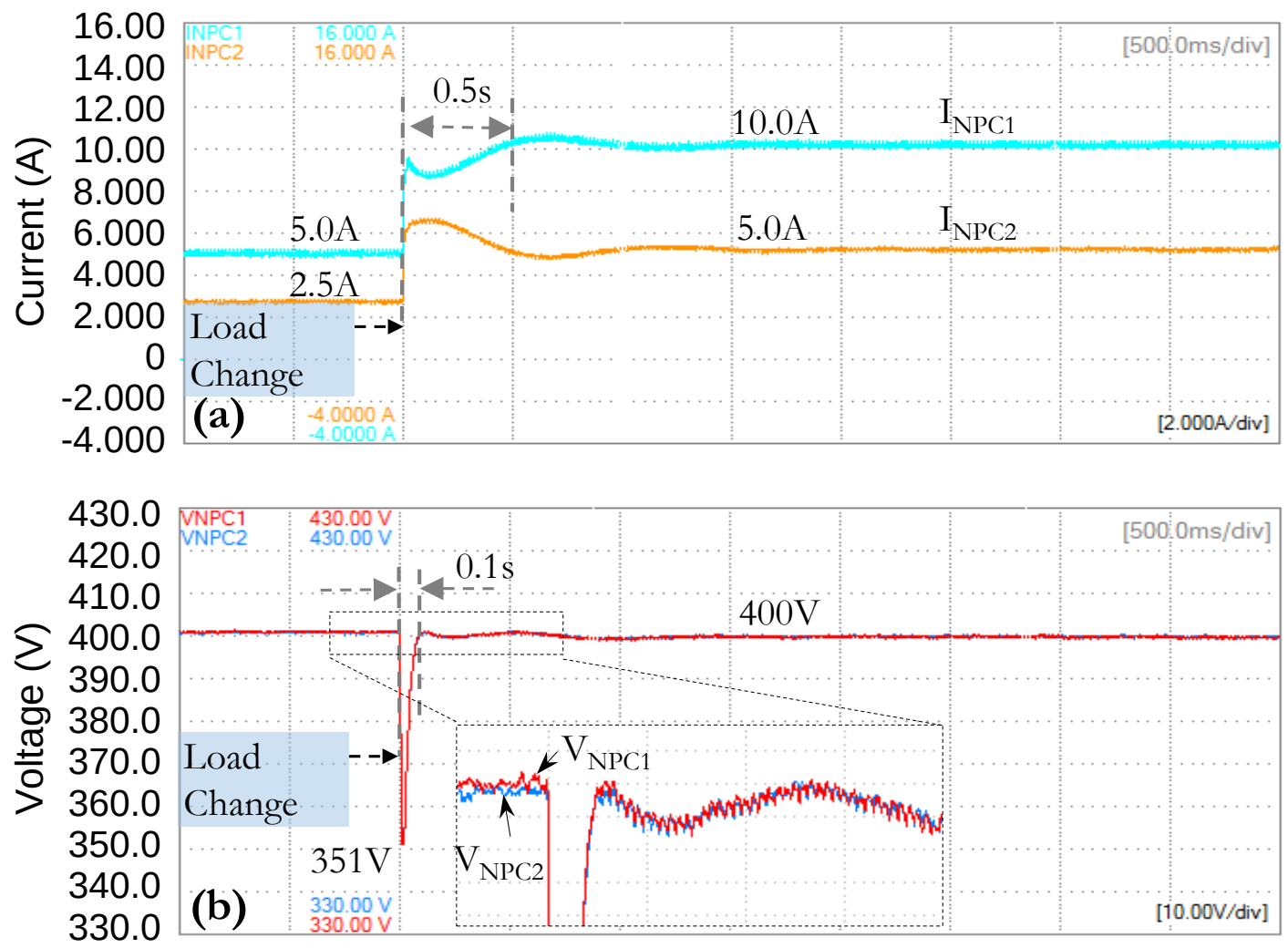

Fig. 12. Nonlinear load change. (a) adaptive current, (b) adaptive voltage. Voltage: 10V/div. Current: 2A/div. Horizontal axis: 500ms/div. Zoomed voltage area (Voltage: 1V/div, Horizontal axis: 200ms/div).

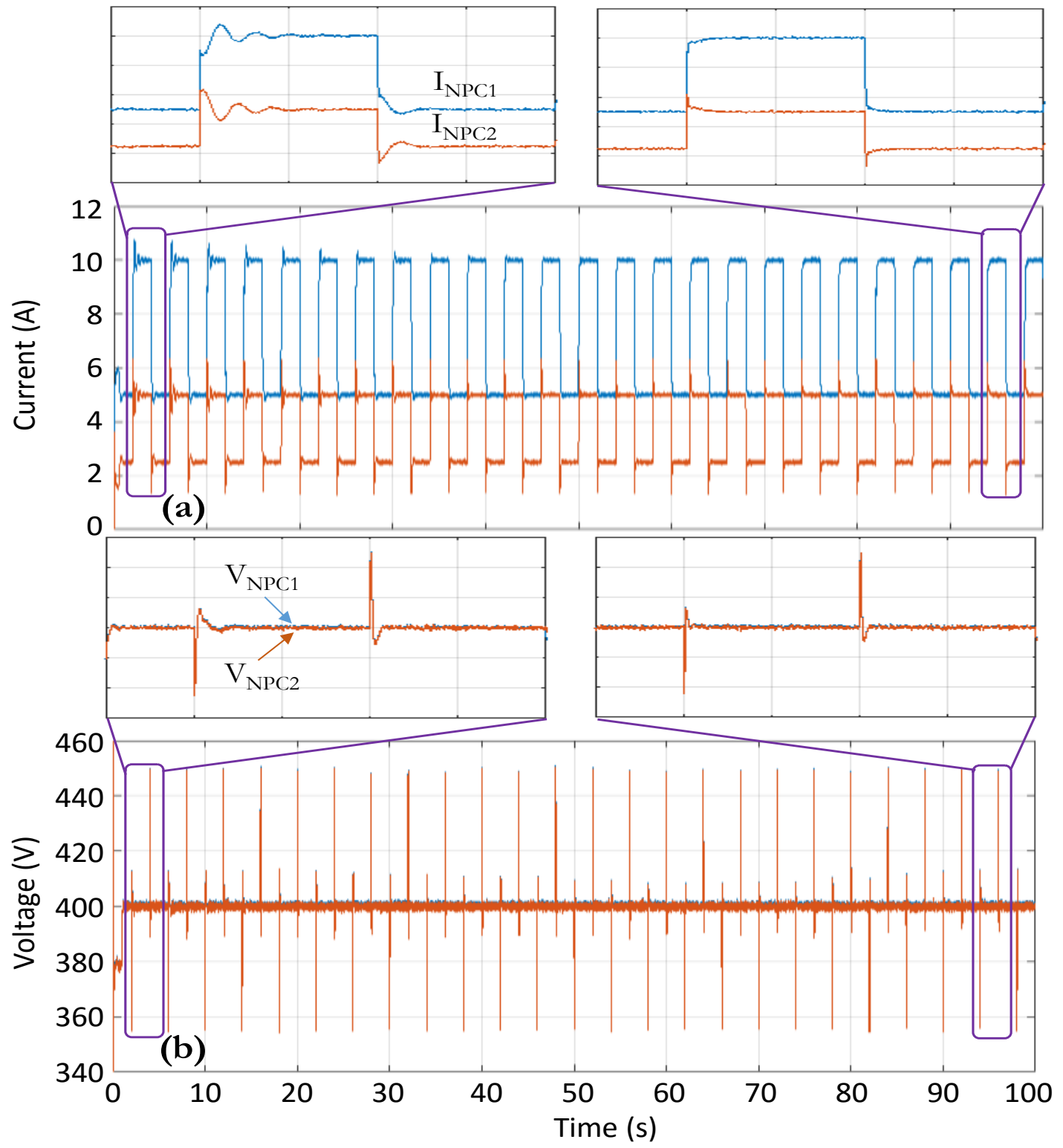

Fig. 13. Adaptive control performance with frequent load changes. (a) current sharing, (b) bus voltage. Zoomed areas (Current: 2A/div, Voltage: 20V/div Horizontal axis: 1s/div).

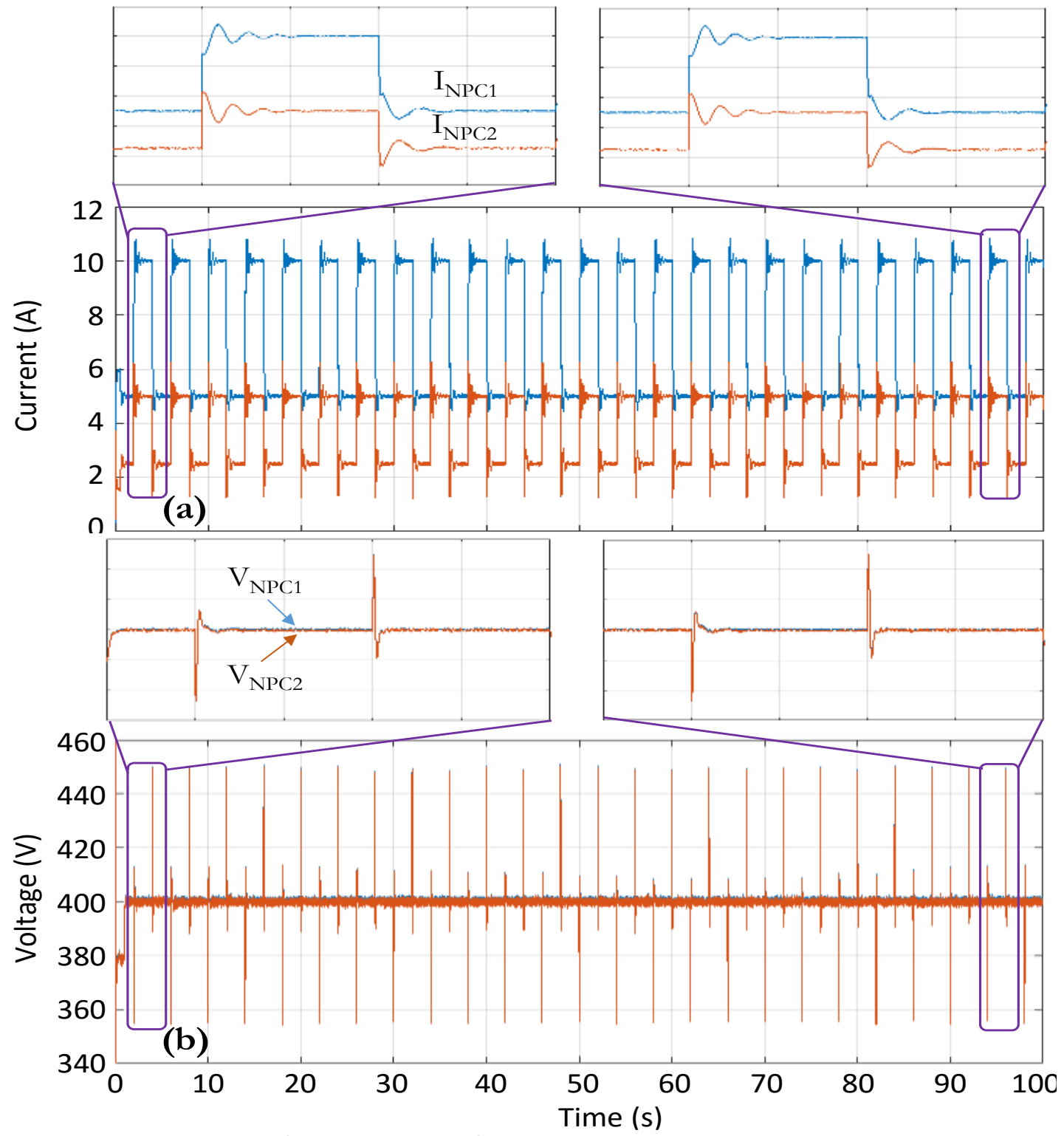

Fig. 14. PI control performance with frequent load changes. (a) current sharing, (b) bus voltage. Zoomed areas (Current: 2A/div, Voltage: 20V/div Horizontal axis: 1s/div).

Consequently, simulation and experimental results show that the proposed method iteratively finds the desired droop resistances to fulfill the system's objectives in current sharing and bus voltage stability in various scenarios.

## 5. Further Demonstration

Additional tests are conducted to further validate the effectiveness of the method in improving the control performance of the system over time. The results are shown in Fig. 13 with arbitrary control parameters selected to present oscillation at the beginning of operation. In this case, the load changes between 3 kW (0.5 pu) and 6 kW (1 pu) every 2 s. Results in Fig. 13 indicate that the current and voltage responses contain more oscillations at the beginning of load changes; however, the oscillations are minimized in the end of the operation (Zoomed areas in Fig. 13). In contrast, Fig. 14 shows that utilizing PI controllers, the system has the same oscillations in the end of operation.

To quantify the improvements in control performance as well as to compare the proposed method with the linear control using PI controllers, the Integral Squared Error (ISE) is utilized for evaluation. The ISE formulas for voltage and current control evaluation are shown in (36) and (37), respectively.

$$ISE_V = \int_0^{\tau} \left( V_{ref} - \frac{V_{NPC1} + V_{NPC2}}{2} \right)^2 dt \tag{36}$$

$$ISE_I = \int_0^{\tau} \left[ \left( I_{ref1} - I_{NPC1} \right)^2 + \left( I_{ref2} - I_{NPC2} \right)^2 \right] dt, \tag{37}$$

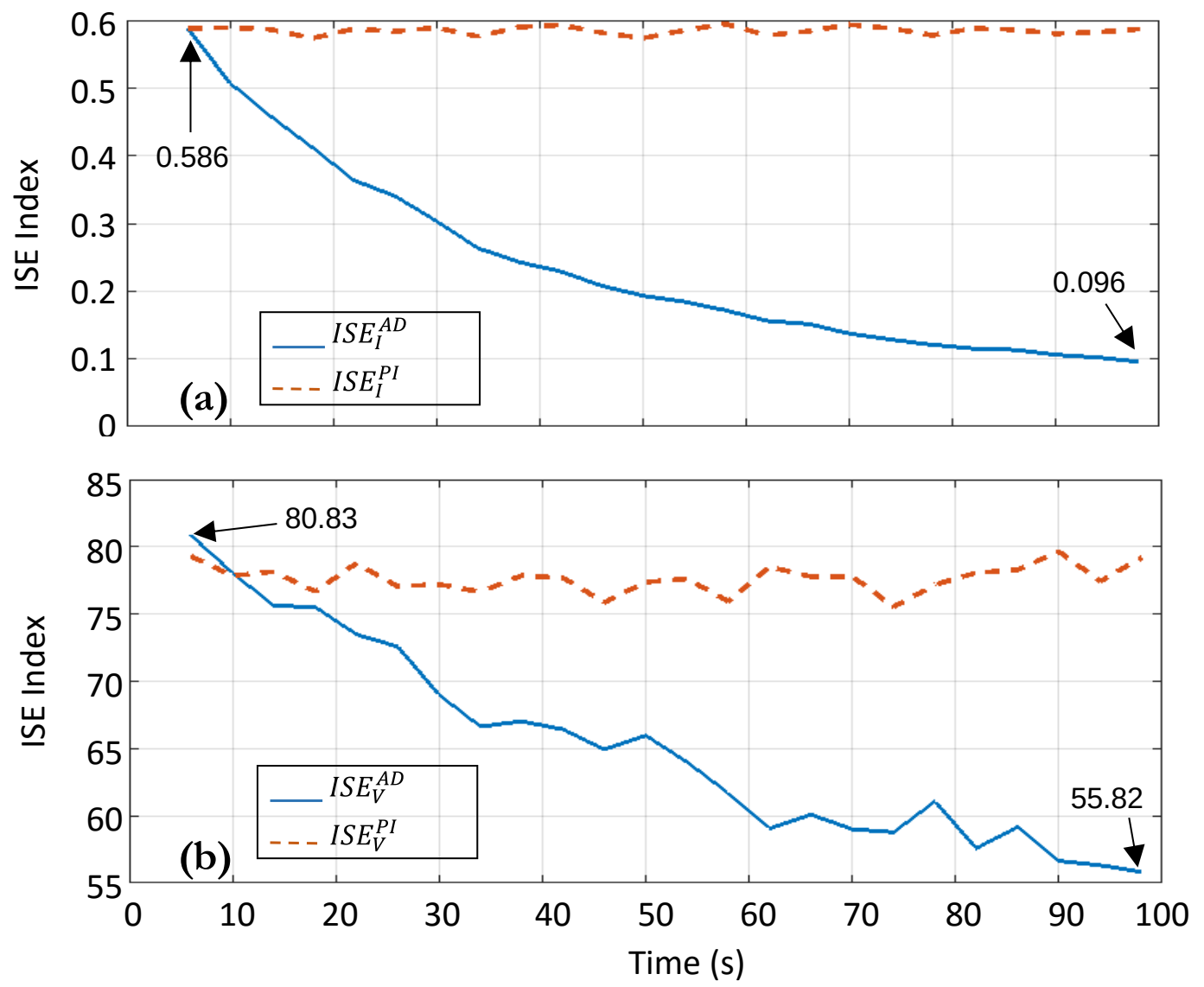

Fig. 15. ISE performance comparison between adaptive control and PI control. (a) ISE indication of current control, (b) ISE indication of voltage control.

where $\tau$ is the evaluation time selected as $\tau = 4\ s$, $V_{ref} = 400\ V$ is the bus voltage reference, $I_{ref1}$ and $I_{ref2}$ are the current sharing references, and $t$ is time.

The comparison indexes for adaptive control are $ISE_I^{AD}$ and $ISE_V^{AD}$; and PI control are $ISE_I^{PI}$ and $ISE_V^{PI}$. These indexes are shown in Fig. 15. In this comparison, the PI control parameters are selected arbitrarily to have similar ISE performance at the beginning of the test. As shown in Fig. 15, the PI control performance indexes $ISE_I^{PI}$ and $ISE_V^{PI}$ vary around their initial values with no indication of improvement in the end of operation. The adaptive control performance indexes $ISE_I^{AD}$ and $ISE_V^{AD}$ demonstrated the improvement over time by the reduction in the ISE values. Specifically, the $ISE_I^{AD}$ decreases to 0.096, which is 16% of its initial value (0.586); and the $ISE_V^{AD}$ decreases to 55.82, which is 69% of its initial value (80.83). Therefore, the proposed adaptive control method iteratively identifies the optimal control parameters to satisfy the voltage and current control objectives of DC microgrids.

## 6. Conclusion

This paper proposes a robust adaptive droop control method for DC microgrids to adjust the droop characteristics to satisfy both power sharing and DC bus voltage stability criteria. We have built comprehensive LTV models to represent the relationship between the droop parameters and the system output in DC microgrids. The closed-loop model reference adaptive control (CRM) is then reviewed, selected, and improved to be applicable to the problem. Stability proof of the proposed improved CRM method based on Lyapunov has been provided in the Appendix. In a distributed consensus framework, simulation and experimental results have validated the capability of the proposed control method to regulate the current sharing and to stabilize the DC bus-voltage of a 400 $V_{dc}$ microgrid system. Additional comparative study using the optimal control index ISE has been conducted to verify the effectiveness of the proposed adaptive algorithm over a linear control method. Future research would consider the communication delays and their effects on the distributed adaptive control algorithm. Additionally, higher order LTV models can be used to compare with the simplified first order adaptive system.

## Appendix

### Stability proof of the adaptation gain-scheduling method

Consider the following Lyapunov function:

$$V = \frac{\tilde{\bar{\theta}}^T \tilde{\bar{\theta}}}{2\gamma_k}|b^*| + \frac{\tilde{b}^2}{2\gamma_k}, \tag{38}$$

where $\tilde{\bar{\theta}} = \bar{\theta} - \bar{\theta}^*$, $\tilde{b} = b - b^*$, and $b^*$, $\bar{\theta}^*$ are the expected values of $\bar{\theta}$ and $b$, respectively. The Lyapunov function has the following derivative:

$$\dot{V} = \frac{2\tilde{\bar{\theta}}^T \dot{\tilde{\bar{\theta}}}}{2\gamma_k}|b^*| + \frac{2\tilde{b}\dot{\tilde{b}}}{2\gamma_k}. \tag{39}$$

The generalization of the control parameter update with parameter projection is

$$\begin{aligned}\dot{\bar{\theta}} &= -\gamma_k \nabla J_{\bar{\theta}} + (1 - sgn|g_{\bar{\theta}}|) \\ &\times \max\left(0, sgn((-\gamma_k sgn(b)\phi_n \varepsilon)^T \bar{\theta})\right)\left(\gamma_k \frac{\nabla g_{\bar{\theta}} \nabla g_{\bar{\theta}}^T}{\nabla g_{\bar{\theta}}^T \gamma \nabla g_{\bar{\theta}}}\right)\gamma_k \nabla J_{\bar{\theta}} \\ \dot{b} &= -\gamma_k \nabla J_b + (1 - sgn|g_b|)\max(0, sgn(\gamma_k(-\bar{\theta}^T \phi_n \\ &\quad + u_n)\varepsilon b))\gamma_k \nabla J_b,\end{aligned} \tag{40}$$

where $J(b, \bar{\theta}) = \frac{m^2\varepsilon^2}{2}$. Substituting (40) into (39), one obtains:

$$\dot{V} = \dot{V}_0 + \dot{V}_{pr}, \tag{41}$$

where

$$\begin{aligned}\dot{V}_0 &= -b^* \varepsilon \tilde{\bar{\theta}}^T \phi_n + \tilde{b}(-\bar{\theta}^T \phi_n + u_n)\varepsilon \\ \dot{V}_{pr} &= (1 - sgn|g_{\bar{\theta}}|) \\ &\times \max\left(0, sgn((-\gamma_k sgn(b)\phi_n \varepsilon)^T \bar{\theta})\right)\tilde{\bar{\theta}}^T \left(\gamma_k \frac{\nabla g_{\bar{\theta}} \nabla g_{\bar{\theta}}^T}{\nabla g_{\bar{\theta}}^T \gamma \nabla g_{\bar{\theta}}}\right)\gamma_k \nabla J_{\bar{\theta}} \\ &+ (1 - sgn|g_b|)\max(0, sgn(\gamma_k(-\bar{\theta}^T \phi_n + u_n)\varepsilon b))\tilde{b}\gamma_k \nabla J_b.\end{aligned} \tag{42}$$

Because

$$\begin{aligned}\varepsilon m^2 &= b^*\left(-\bar{\theta}^{*T}\phi_n + u_n\right) - b(-\bar{\theta}^T \phi_n + u_n) \\ &= -b^* \bar{\theta}^{*T} \phi_n + b\bar{\theta}^T \phi_n - \tilde{b} u_n + b^* \bar{\theta}^T \phi_n - b^* \bar{\theta}^T \phi_n \\ &= b^* \tilde{\bar{\theta}}^T \phi_n - \tilde{b}(-\bar{\theta}^T \phi_n + u_n),\end{aligned} \tag{43}$$

then

$$\dot{V}_0 = -b^* \varepsilon \tilde{\bar{\theta}}^T \phi_n + \left(b^* \tilde{\bar{\theta}}^T \phi_n - \varepsilon m^2\right)\varepsilon = -\varepsilon^2 m^2 \le 0. \tag{44}$$

To prove $\dot{V}_{pr} \le 0$, first we define that $\Omega_{0\bar{\theta}}$ and $\Omega_{0b}$ are the boundaries of convex sets $M_{\bar{\theta}}$ and $M_b$ and the interiors of $M_{\bar{\theta}}$ and $M_b$ are $\Omega_{\bar{\theta}}$ and $\Omega_b$. Following cases are considered for $\dot{V}_{pr}$:

Case 1: If $\bar{\theta} \in \Omega_{0\bar{\theta}}$ and $b \in \Omega_{0b}$. The function $\dot{V}_{pr} = 0$.

Case 2: If $\bar{\theta} \in \Omega_{0\bar{\theta}}$ and $b \notin \Omega_{0b}$ the following achieves:

$$\begin{aligned}\dot{V}_{pr} &= \tilde{\bar{\theta}}^T \left(\gamma_k \frac{\nabla g_{\bar{\theta}} \nabla g_{\bar{\theta}}^T}{\nabla g_{\bar{\theta}}^T \gamma_k \nabla g_{\bar{\theta}}}\right)\gamma_k \nabla J_{\bar{\theta}} \\ &= \frac{1}{\nabla g_{\bar{\theta}}^T \gamma_k \nabla g_{\bar{\theta}}}(\tilde{\bar{\theta}}^T \nabla g_{\bar{\theta}})((\gamma_k \nabla J_{\bar{\theta}})^T \nabla g_{\bar{\theta}}).\end{aligned} \tag{45}$$

The expected value $\bar{\theta}^* \in M_{\bar{\theta}}$, and value $\bar{\theta} \in \Omega_{0\bar{\theta}}$. Therefore, vector $\tilde{\bar{\theta}}^T$ is in the same direction to vector $\nabla g_{\bar{\theta}}$, which results in $\tilde{\bar{\theta}}^T \nabla g_{\bar{\theta}} \ge 0$. Moreover, the convex property of $M_{\bar{\theta}}$ makes $((\gamma_k \nabla J_{\bar{\theta}})^T \nabla g_{\bar{\theta}} < 0$. Hence, $\dot{V}_{pr} \le 0$ in this case.

Case 3: If $b \in \Omega_{0b}$ and $\bar{\theta} \notin \Omega_{0\bar{\theta}}$, and one attains:

$$\dot{V}_{pr} = \gamma_k \tilde{b} \nabla J_b. \tag{46}$$

The parameter $\tilde{b}$ develops in the opposite direction to $\nabla J_b$, which results in $\dot{V}_{pr} = \gamma_k \tilde{b} \nabla J_b \le 0$.

Consequently, $\dot{V}_{pr} \le 0$ in all cases, which results in $\dot{V} = \dot{V}_0 + \dot{V}_{pr} \le 0$. Moreover, $V$ is positive definite. Therefore, the system is stable in the sense of Lyapunov.